\documentclass[12pt,fleqn]{article}
\usepackage{mathrsfs}
\usepackage{amsfonts}
\usepackage{amssymb}
\usepackage{latexsym,amsmath}
\usepackage{graphicx}
\date{}
\begin{document}
\title{The Q-generating function for graphs with application}
\author{Shu-Yu Cui$^a$, Gui-Xian Tian$^b$\footnote{Corresponding author.  E-mail: gxtian@zjnu.cn or guixiantian@gmail.com (G.-X.
Tian)}\\
{\small{\it $^a$ Xingzhi College, Zhejiang Normal University,
Jinhua, Zhejiang, 321004,
P.R. China}}\\
{\small{\it $^b$College of Mathematics, Physics and Information Engineering,}}\\
{\small{\it Zhejiang Normal University, Jinhua, Zhejiang, 321004,
P.R. China}}}\maketitle

\begin{abstract}

For a simple connected graph $G$, the $Q$-generating function of the
numbers $N_k$ of semi-edge walks of length $k$ in $G$ is defined by
$W_Q(t)=\sum\nolimits_{k = 0}^\infty  {N_k t^k }$. This paper
reveals that the $Q$-generating function $W_Q(t)$ may be expressed
in terms of the $Q$-polynomials of the graph $G$ and its complement
$\overline{G}$. Using this result, we study some $Q$-spectral
properties of graphs and compute the $Q$-polynomials for some graphs
obtained by the use of some operation on graphs, such as the
complement graph of a regular graph, the join of two graphs, the
(edge)corona of two graphs and so forth. As another application of
the $Q$-generating function $W_Q(t)$, we also give a combinatorial
interpretation of the $Q$-coronal of $G$, which is defined to be the
sum of the entries of the matrix $(\lambda I_n-Q(G))^{-1}$. This
result may be used to obtain the many alternative calculations of
the $Q$-polynomials of the (edge)corona of two graphs. Further, we
also compute the $Q$-coronals of the join of two graphs and the
complete multipartite graphs.

\emph{AMS classification:} 05C50 05C90

\emph{Keywords:} Signless Laplacian matrix; $Q$-polynomial;
$Q$-Spectrum; $Q$-generating function; $Q$-Coronal; Semi-edge walk

\end{abstract}

\section*{1. Introduction}

\indent\indent Throughout this paper, we consider only simple
connected graphs. Let $G=(V,E)$ be a graph with vertex set
$V=\{v_1,v_2,\ldots,v_n\}$. Two vertices $v_i$ and $v_j$ of $G$ are
called \emph{adjacent}, denoted by $v_i\sim v_j$, if they are
connected by an edge. The \emph{adjacency matrix} $A(G)$ of $G$ is a
square matrix of order $n$, whose entry $a_{i,j}$ is defined as
follows: $a_{i,j}=1$ if $v_i\sim v_j$, $0$ otherwise. Let $D(G)$ be
the diagonal degree matrix of $G$. The matrix $Q(G)=D(G)+A(G)$ is
called the \emph{signless Laplacian matrix} of $G$. The
\emph{$Q$-spectrum} of $G$ is defined to
\[
S(G) = (q_1 (G),q_2 (G), \ldots ,q_n (G)),
\]
where $q_1 (G)\geq q_2 (G)\geq\cdots\geq q_n (G)$ are the
eigenvalues of $Q(G)$. They also are the roots of the
\emph{$Q$-polynomial} $f_Q(\lambda)=\det(\lambda I_n-Q(G))$ of $G$.
Denote $Q$-polynomial of the complement graph $\overline{G}$ of $G$
by $f_{\overline{Q}}(\lambda)=\det(\lambda I_n-Q(\overline{G}))$.
For more review about the $Q$-spectrum and $Q$-polynomial of $G$,
readers may refer to
\cite{Aouchiche2010,Cvetkovic2008,Cvetkovic2007,CvetkovicI,CvetkovicII,CvetkovicIII,Cvetkovic1980,Oliveira2010}
and the references therein.

Let $G$ be a simple connected graph and $A(G)$ be its adjacency
matrix. A \emph{walk} (of length $k$) in $G$ is an alternating
sequence $v_1,e_1,v_2,e_2,\ldots,v_k, e_k, v_{k+1}$ of vertices
$v_1, v_2,\ldots, v_{k+1}$ and edges $e_1, e_2,\ldots, e_k$ such
that for any $i=1,2,\ldots, k$ the vertices $v_i$ and $v_{i+1}$ are
distinct end-vertices of the edge $e_i$. It is well
known\cite{Cvetkovic1980} that the $(i,j)$-entry of the matrix
$A(G)^k$ equals the number of walks of length $k$ starting at vertex
$v_i$ and terminating at vertex $v_j$. Let $N_k$ denote the total
number of all walks of length $k$ in $G$. $H_G(t)=\sum\nolimits_{k =
0}^\infty  {N_k t^k }$ is called the \emph{generating function} of
the numbers $N_k$ of all walks of length $k$ in $G$. In
\cite{Cvetkovic1980}, the generating function $H_G(t)$ is expressed
in terms of the characteristic polynomials of the graph $G$ and its
complement $\overline{G}$, and many spectral properties are
obtained. For example, the characteristic polynomials of some graphs
is computed by employing the generating function $H_G(t)$ in
\cite{Cvetkovic1980}.

For a simple connected graph $G$, let $Q(G)$ be its signless
Laplacian matrix. Similarly, a \emph{semi-edge walk} (of length $k$)
\cite{Cvetkovic2007} in an (undirected) graph $G$ is an alternating
sequence $v_1,e_1,v_2,e_2,\ldots,v_k, e_k, v_{k+1}$ of vertices
$v_1, v_2,\ldots, v_{k+1}$ and edges $e_1, e_2,\ldots, e_k$ such
that for any $i=1,2,\ldots, k$ the vertices $v_i$ and $v_{i+1}$ are
end-vertices (not necessarily distinct) of the edge $e_i$. It is
proved \cite{Cvetkovic2007} that the $(i,j)$-entry of the matrix
$Q(G)^k$ equals the number of semi-edge walks of length $k$ starting
at vertex $v_i$ and terminating at vertex $v_j$.

The \emph{$Q$-generating function} of the numbers $N_k$ of semi-edge
walks of length $k$ in $G$ is defined to $W_Q(t)=\sum\nolimits_{k =
0}^\infty  {N_k t^k }$, where $N_k$ denotes the total number of
semi-edge walks of length $k$ in $G$. The following problem seems
interesting:

\emph{Study the $Q$-generating function $W_Q(t)$ for the numbers
$N_k$ of semi-edge walks of length $k$ in $G$ and compute the
$Q$-polynomials of some graphs by employing the $Q$-generating
function $W_Q(t)$.}

This paper reveals that the $Q$-generating function $W_Q(t)$ may be
expressed in terms of the $Q$-polynomials of the graph $G$ and its
complement $\overline{G}$. Using this result, we obtain some
$Q$-spectral properties of graphs and compute the $Q$-polynomials
for some graphs obtained by the use of some operation on graphs,
such as the complement of a graph, the join of two graphs, the
(edge)corona of two graphs and so forth.

As another application of the $Q$-generating function $W_Q(t)$, we
also give a combinatorial interpretation of the $Q$-coronal of a
graph $G$, which is defined to be the sum of the entries of the
matrix $(\lambda I_n-Q(G))^{-1}$. This result may be used to obtain
the many alternative calculations of the $Q$-polynomials of the
(edge)corona of two graphs. Further, we also compute the
$Q$-coronals of the join of two graphs and complete multipartite
graphs.

\section*{2. The $Q$-generating functions and $Q$-polynomials of graphs}

\indent\indent For a simple connected graph $G$, the following
Proposition 2.1 reveals that the $Q$-generating function $W_Q(t)$
may be expressed in terms of the $Q$-polynomials
of the graph $G$ and its complement $\overline{G}$.\\
\\
 \textbf{Proposition 2.1.} {\it Let $G$ be a simple connected graph
on $n$ vertices. Then
\[
W_Q (t) = \frac{1}{t}\left( {( - 1)^n \frac{{f_{ \overline{Q} }
\left( {n - 2 - \frac{1}{t}} \right)}}{{f_Q \left( {\frac{1}{t}}
\right)}} - 1} \right).
\]}\\
\textbf{Proof.} The proof is totally similar to Theorem 1.11 in
\cite{Cvetkovic1980}. Let $B$ be a nonsingular $n$-by-$n$ square
matrix and $J$ be a square matrix all entries of which are equal to
$1$. Then, for arbitrary number $x$,
\begin{equation}\label{1}
\det (B + xJ) = \det B + x\;{\rm{sum}}({\rm{adj}}B),
\end{equation}
where $\text{sum}(K)$ denotes the sum of all entries of a matrix $K$
and $\text{adj}K$ denotes its adjoint matrix.

Now, from Theorem 4.1 in \cite{Cvetkovic2007},  one gets
$N_k=\text{sum}(Q^k)$. Noting that
\[
\sum\limits_{k = 0}^\infty  {Q^k t^k }  = (I - tQ)^{ - 1}  =
\frac{{\text{adj}(I - tQ)}}{{\det (I - tQ)}}\quad\quad\quad\left(|t|
< \frac{1}{{q_1 }}\right).
\]
Thus we obtain
\[
W_Q (t) = \sum\limits_{k = 0}^\infty  {N_k t^k }  = \sum\limits_{k =
0}^\infty  {\text{sum}(Q^k )t^k }  = \frac{{\text{sum}(\text{adj}(I
- tQ))}}{{\det (I - tQ)}}.
\]
With $B=I-tQ$, $x=t$, the formula (\ref{1}) yields
\begin{align*}
\text{sum}(\text{adj}(I - tQ))& = \frac{1}{t}\left( {\det (I - tQ + tJ) - \det (I - tQ)} \right) \\
&= \frac{1}{t}\left( {\det ((1 - (n - 2)t)I + t\overline Q ) - \det
(I - tQ)} \right),
\end{align*}
where $\overline{Q}=(n-2)I+J-Q$ is the signless Laplacian matrix of
the complement $\overline{G}$ of $G$. Hence,
\begin{align*}
W_Q (t) &= \frac{1}{t}\left( {\frac{{\det ((1 - (n - 2)t)I +
t\overline Q )}}{{\det (I - tQ)}} - 1} \right) \\
&= \frac{1}{t}\left( {( - 1)^n \frac{{\det \left( {(n - 2 -
\frac{1}{t})I - \overline Q } \right)}}{{\det \left( {\frac{1}{t}I -
Q} \right)}} - 1} \right) \\
&= \frac{1}{t}\left( {( - 1)^n \frac{{f_{\overline Q } \left( {n - 2
- \frac{1}{t}} \right)}}{{f_Q \left( {\frac{1}{t}} \right)}} - 1}
\right).
\end{align*}
This completes the proof of Proposition 2.1. $\Box$\\
\\
\textbf{Theorem 2.2.} {\it Let $G$ be an $r$-regular graph on $n$
vertices and $\overline{G}$ be its complement graph. Then
\[
f_{\overline Q } (\lambda ) = ( - 1)^n \left( {1 + \frac{n}{{n - 2 -
2r - \lambda }}} \right)f_Q (n - 2 - \lambda ).
\]
Moreover, if the signless Laplacian spectrum of $G$ contains
$2r,q_2,\ldots,q_n$, then the signless Laplacian spectrum of
$\overline{G}$ contains $2(n-r-1),n-2-q_2,\ldots,n-2-q_n$.}\\
\\
\textbf{Proof.} It is easy to see that
\[
W_Q (t) = \sum\limits_{k = 0}^\infty  {N_k t^k }  = \sum\limits_{k =
0}^\infty  {n(2r)^k t^k }  = \frac{n}{{1 - 2rt}}\quad \quad\left(
{|t| < \frac{1}{{2r}}} \right),
\]
whenever $G$ is an $r$-regular graph on $n$ vertices. From
Proposition 2.1, one has
\begin{align}\label{2}
\frac{1}{t}\left( {( - 1)^n \frac{{f_{\overline Q } \left( {n - 2 -
\frac{1}{t}} \right)}}{{f_Q \left( {\frac{1}{t}} \right)}} - 1}
\right) = \frac{n}{{1 - 2rt}}.
\end{align}
With $\lambda=(n-2)-\frac{1}{t}$, that is,
$\frac{1}{t}=(n-2)-\lambda$ in (\ref{2}), we obtain the required
result. Moreover, if the signless Laplacian spectrum of $G$ contains
$2r,q_2,\ldots,q_n$, then it is easy to see that the signless
Laplacian spectrum of $\overline{G}$ contains
$2(n-r-1),n-2-q_2,\ldots,n-2-q_n$. $\Box$\\

Let $G_1$ and $G_2$ be two graphs with disjoint vertex sets $V(G_1)$
and $V(G_2)$, and edge sets $E(G_1)$ and $E(G_2)$, respectively. The
\emph{join} $G_1\vee G_2$ of $G_1$ and $G_2$ is the graph union
$G_1\cup G_2$ together with all the edges joining $V(G_1)$ and
$V(G_2)$.\\
\\
\textbf{Theorem 2.3.} {\it Let $G_1$ and $G_2$ be two simple
connected graphs with $n_1$ and $n_2$ vertices, respectively. Then
\begin{align*}
f_{Q(G_1  \vee G_2 )} = &\; ( - 1)^{n_2 } f_{Q_1 } (\lambda  - n_2
)f_{\overline {Q_2 } } (n_1  + n_2  - \lambda  - 2)\\
& + ( - 1)^{n_1 } f_{Q_2 } (\lambda  - n_1 )f_{\overline {Q_1 } }
(n_1 + n_2  - \lambda  - 2)\\
& - ( - 1)^{n_1  + n_2 } f_{\overline {Q_1 } } (n_1  + n_2  -
\lambda - 2)f_{\overline {Q_2 } } (n_1  + n_2  - \lambda  - 2),
\end{align*}
} where $Q_1,Q_2$ are the signless Laplacian matrices of $G_1$ and
$G_2$, respectively.\\
\\
\textbf{Proof.} Clearly, $W_{Q(G_1\oplus
G_2)}(t)=W_{Q_1}(t)+W_{Q_2}(t)$, where $G_1\oplus G_2$ denotes the
direct sum of $G_1$ and $G_2$. Proposition 2.1 implies that
\begin{align}\label{3}
&\frac{1}{t}\left( {( - 1)^{n_1  + n_2 }  \frac{{f_{\overline Q (G_1
\oplus G_2 )} \left( {n_1  + n_2  - 2 - \frac{1}{t}}
\right)}}{{f_{Q_1 } \left( {\frac{1}{t}} \right)f_{Q_2 } \left(
{\frac{1}{t}} \right)}}
- 1} \right) \nonumber\\
&\quad\quad\quad\quad\quad\quad= \sum\limits_{i = 1}^2
{\frac{1}{t}\left( {( - 1)^{n_i } \frac{{f_{\overline {Q_i } }
\left( {n_i  - 2 - \frac{1}{t}} \right)}}{{f_{Q_i } \left(
{\frac{1}{t}} \right)}} - 1} \right)}.
\end{align}
Note that $\overline{G_1\oplus G_2}=\overline{G_1}\vee
\overline{G_2}$. Setting ${n_1  + n_2  - 2 - \frac{1}{t}}=\lambda$
and substituting $\overline{G_1}$, $\overline{G_2}$ for $G_1$, $G_2$
in (\ref{3}), we
obtain the required result. $\Box$\\
\\
\textbf{Corollary 2.4}\cite{Freitas2010}. {\it Let $G_i$ $(i=1,2)$
be a regular graph of degree $r_i$ with $n_i$ vertices. Then
\[
f_{Q(G_1  \vee G_2 )}  = \left( {1 - \frac{{n_1 n_2 }}{{(\lambda  -
n_1  - 2r_2 )(\lambda  - n_2  - 2r_1 )}}} \right) f_{Q_1 } (\lambda
- n_2 )f_{Q_2 } (\lambda  - n_1 ).
\]
}\\
\textbf{Proof.} This is an immediate consequence
of Theorems 2.2 and 2.3, omitted. $\Box$\\

Let $Q(G)=(q_{ij})_{n\times n}$ be the signless Laplacian of a
simple graph $G$. Assume that $x_1,x_2,\ldots,x_n$ are mutually
orthogonal normalized eigenvectors of $Q(G)$ associated to
eigenvalues $q_1,q_2,\ldots,q_n$, respectively. Also let
$\Lambda=\text{diag}(q_1,q_2,\ldots,q_n)$ and
$P=(x_1,x_2,\ldots,x_n)=(x_{ij})_{n\times n}$. Then $Q(G)=P\Lambda
P^T$, which implies that the number $N_k$ of all semi-edge walks of
length $k$ in $G$ equals
\[
N_k  = \sum\limits_{i,j} {N_k (i,j)}  = \sum\limits_{i,j}
{q_{ij}^{(k)} }  = \sum\limits_{l = 1}^n {\left( {\sum\limits_{i =
1}^n {x_{il} } }\right )^2 q_l^k }.
\]
Thus we arrive at:\\
\\
\textbf{Theorem 2.5.} {\it The total number $N_k$ of semi-edge walks
of length $k$ in $G$ equals
\[
N_k  = \sum\limits_{l = 1}^n {\gamma _l q_l^k } \;\;\;\; (k = 0,1,2,
\ldots ),
\]
where $\gamma _l  = \left( {\sum\limits_{i = 1}^n {x_{il} } }
\right)^2.$}\\

It is clear to see that $N_k=n(2r)^k$ whenever $G$ is an $r$-regular
graph with $n$ vertices. In this case, the signless Laplacian
spectral radius $q_1$ of $G$ is equal to
\[
q_1  = \sqrt[k]{{\frac{{N_k }}{n}}} = 2r.
\]
In general case, we have the following Theorem 2.6, which is
analogous to an existing result related to the adjacency spectrum
(see Theorem 1.12 in \cite{Cvetkovic1980}).
\\
\\
\textbf{Theorem 2.6.} {\it
\[
q_1  = \mathop {\lim }\limits_{k \to \infty } \sqrt[k]{{\frac{{N_k
}}{n}}}=\mathop {\lim }\limits_{k \to \infty } \sqrt[k]{{N_k }}.
\]
}\\
\\
\textbf{Proof.} Firstly, it is easy to see that
\[
q_1 \sqrt[k]{{\frac{{\gamma _1 }}{n}}} \le \sqrt[k]{{\frac{{N_k
}}{n}}} = \sqrt[k]{{\frac{{\sum\nolimits_{l = 1}^n {\gamma _l q_l^k
} }}{n}}} \le q_1 \sqrt[k]{{\frac{{\sum\nolimits_{l = 1}^n {\gamma
_l } }}{n}}}.
\]
The Squeeze Theorem implies that
\[
q_1  = \mathop {\lim }\limits_{k \to \infty } \sqrt[k]{{\frac{{N_k
}}{n}}}.
\]
Note that $\mathop {\lim }\nolimits_{k \to \infty } \sqrt[k]{{n
}}=1$, the required result follows. $\Box$\\

The following statement and its proof is analogous to an existing
result related to the adjacency spectrum (see Theorem 2.5 in
\cite{Cvetkovic1980}).\\
\\
\textbf{Theorem 2.7.} {\it If the Q-spectrum of a graph $G$ contains
a signless Laplacian eigenvalue $q_0$ with multiplicity $s\geq 2$,
then the Q-spectrum of its complementary graph $\overline{G}$
contains a signless Laplacian eigenvalue $n-2-q_0$ with multiplicity
$t$, where $s-1\leq t\leq s+1$.}\\
\\
\textbf{Proof.}  By Theorem 2.5, the $Q$-generating function of the
numbers $N_k$ of semi-edge walks of length $k$ in $G$ is
\begin{align*}
W_Q (t)& = \sum\limits_{k = 0}^\infty  {N_k t^k }  = \sum\limits_{k
= 0}^\infty  {\left( {\sum\limits_{l = 1}^n {\gamma _l q_l^k } }
\right)t^k } \\
& = \sum\limits_{l = 1}^n {\gamma _l \left( {\sum\limits_{k =
0}^\infty  {q_l^k t^k } } \right)} = \sum\limits_{l = 1}^n
{\frac{{\gamma _l }}{{1 - tq_l }}}\quad\quad\quad \left( {|t| <
\frac{1}{{q_1 }}} \right).
\end{align*}
Set
\[
\Phi (u) = ( - 1)^n \frac{{f_{\overline Q } \left( {n - 2 - u}
\right)}}{{f_Q (u)}}.
\]
From Proposition 2.1, one has
\[
\Phi (u) = 1 + \frac{1}{u}W_Q \left( {\frac{1}{u}} \right) = 1 +
\sum\limits_{l = 1}^n {\frac{{\gamma _l }}{{u - q_l }}}  =
\frac{{\varphi _1 (u)}}{{\varphi _2 (u)}},
\]
where $\varphi _1 (u)$, $\varphi _2 (u)$ are polynomials in $u$ and
the roots of $\varphi _2 (u)$ are all simple. Since $q_0$ is a
signless Laplacian eigenvalue of $G$ with multiplicity $s\geq 2$.
Then $f_Q(u)=(u-q_0)^sg(u)$, where $g(q_0)\neq 0$. Therefore,
\[
\Phi (u) = ( - 1)^n \frac{{f_{\overline Q } \left( {n - 2 - u}
\right)}}{{(u - q_0 )^s g(u)}} = \frac{{\varphi _1 (u)}}{{\varphi _2
(u)}},
\]
which implies that $f_{\overline Q } \left( {n - 2 - u} \right)$
must have a factor $(u - q_0 )^t$, $t\geq s-1$ as the roots of
$\varphi _2 (u)$ are all simple.  Thus $f_{\overline Q } (u )$ must
have a factor $(u - (n-2-q_0 ))^t$. Hence the Q-spectrum of the
complementary graph $\overline{G}$ contains a signless Laplacian
eigenvalue $n-2-q_0$ with multiplicity $t\geq s-1$.

Next, we shall prove $t\leq s+1$. Assume that the Q-spectrum of the
complementary graph $\overline{G}$ contains a signless Laplacian
eigenvalue $n-2-q_0$ with multiplicity $t> s+1$. According to the
above statement, $\overline{\overline{G}}=G$ contains a signless
Laplacian eigenvalue $n-2-(n-2-q_0)=q_0$ with multiplicity $r\geq
t-1>s+1-1=s$, contradiction. Hence, $t\leq s+1$.

This completes the proof of Theorem. $\Box$\\

Next, we shall consider another applications of the $Q$-generating
function $W_Q(t)$ of the numbers $N_k$ of semi-edge walks of length
$k$ in $G$. In \cite{Cui2012}, Cui and Tian introduced a new
invariant, the \emph{$Q$-coronal} $\Gamma_Q(\lambda)$ of a graph $G$
of order $n$. It is defined to be the sum of the entries of the
matrix $(\lambda I_n-Q)^{-1}$, where $I_n$ and $Q$ are the identity
matrix of order $n$ and the signless Laplacian matrix of $G$,
respectively. Using this concept, we computed the $Q$-polynomials of
the corona $G_1\circ G_2$ and edge corona $G_1\diamond G_2$ (for
definitions and more details about the corona and edge corona, see
\cite{Barik2007,Cui2012,Cui20122,Hou2010,Liu2013,McLeman2011}) as
follows.
\\
\\
\textbf{Theorem 2.8}\cite{Cui2012}. {\it Let $G_1$ and $G_2$ be two
graphs on $n_1$ and $n_2$ vertices, respectively. Also let $\Gamma
_{Q_2} (\lambda )$ be the $Q_2$-coronal of $G_2$ and $G=G_1\circ
G_2$. Then the $Q$-polynomial of $G$ is
\[
f_Q (\lambda ) = (f_{Q_2 } (\lambda  - 1))^{n_1 } f_{Q_1 } (\lambda
- n_2 - \Gamma _{Q_2 } (\lambda  - 1)).
\]
}\\
\textbf{Theorem 2.9}\cite{Cui2012}. {\it Let $G_1$ be an
$r_1$-regular graph with $n_1$ vertices, $m_1$ edges and $G_2$ be
any graph with $n_2$ vertices, $m_2$ edges. Also let $\Gamma _{Q_2}
(\lambda )$ be the $Q_2$-coronal of $G_2$ and $G=G_1\diamond G_2$.
If $\lambda$ is not a pole of $\Gamma _{Q_2} (\lambda-2)$, then the
$Q$-polynomial of $G$ is
\[
f_Q (\lambda ) = (f_{Q_2 } (\lambda  -2))^{m_1 } f_{Q_1 }
\left(\frac{\lambda -r_1 n_2}{1+ \Gamma _{Q_2 } (\lambda  - 2)}
\right)(1+ \Gamma _{Q_2 } (\lambda  - 2))^{n_1}.
\]}

It is well known that it is difficult for us to compute the inverse
of matrices, especially high order matrices, which results in a
difficulty when we need to compute the $Q$-coronal
$\Gamma_Q(\lambda)$ in Theorems 2.8 and 2.9. In \cite{Cui2012}, we
computed the $Q$-coronal of some special graphs and gave the
Q-polynomials of their (edge)coronae.

Next, we shall give a combinatorial interpretation of the
$Q$-coronal of a graph $G$ of order $n$, which is used to obtain the
many alternative calculations of the $Q$-polynomials of the corona
$G_1\circ G_2$ and edge corona $G_1\diamond G_2$ for any graphs
$G_1$ and $G_2$.\\
\\
\textbf{Proposition 2.10.} {\it Let $G$ be a simple connected graph
of order $n$. Then it's $Q$-coronal equals
\[
\Gamma _Q (\lambda ) =  - 1 + ( - 1)^n \frac{{f_{\overline Q }
\left( {n - 2 - \lambda } \right)}}{{f_Q (\lambda )}}.
\]}\\
\textbf{Proof.} Let $Q$ be the signless Laplacian matrix of $G$ and
$\textbf{1}_n$ denote the length-$n$ column vector, whose each
element equals $1$. By a simple calculation,
\begin{align}\label{4}
\Gamma _Q (\lambda )= \textbf{1}_n^T (\lambda I_n  - Q)^{ - 1}
\textbf{1}_n & = \lambda ^{ - 1} \textbf{1}_n^T (I_n  - \lambda ^{ -
1} Q)^{ - 1} \textbf{1}_n \nonumber\\& = \frac{1}{\lambda
}\textbf{1}_n^T \left( {\sum\limits_{k = 0}^\infty  {Q^k \left(
{\frac{1}{\lambda }} \right)^k } } \right)\textbf{1}_n \nonumber
\\& = \frac{1}{\lambda }\sum\limits_{k = 0}^\infty  {\left( {\textbf{1}_n^T
Q^k \textbf{1}_n } \right)\left( {\frac{1}{\lambda }} \right)^k }.
\end{align}
Since the sum $\textbf{1}_n^T Q^k \textbf{1}_n $ of all elements of
$Q^k$ is the total number $N_k$ of all semi-edge walks of length $k$
in $G$. Then the equality (\ref{4}) becomes
\[
\Gamma _Q (\lambda ) = \frac{1}{\lambda }\sum\limits_{k = 0}^\infty
{N_k \left( {\frac{1}{\lambda }} \right)^k }  = \frac{1}{\lambda
}W_Q \left( {\frac{1}{\lambda }} \right).
\]
From Proposition 2.1, the required result follows. $\Box$\\

Now applying Proposition 2.10, Theorems 2.8 and 2.9 may be rewritten
as the following Theorems 2.11 and 2.12, respectively.
\\
\\
\textbf{Theorem 2.11.} {\it Let $G_1$ and $G_2$ be two graphs on
$n_1$ and $n_2$ vertices, respectively. Also let $G=G_1\circ G_2$.
Then the $Q$-polynomial of $G$ is
\[
f_Q (\lambda ) = (f_{Q_2 } (\lambda  - 1))^{n_1 } f_{Q_1
}\left(\lambda - n_2 + 1 -( - 1)^{n_2} \frac{{f_{\overline Q_2 }
\left( {n - \lambda -1} \right)}}{{f_{Q_2} (\lambda-1 )}}\right).
\]
}\\
\textbf{Theorem 2.12.} {\it Let $G_1$ be an $r_1$-regular graph with
$n_1$ vertices, $m_1$ edges and $G_2$ be any graph with $n_2$
vertices, $m_2$ edges. Also let $G=G_1\diamond G_2$. Then the
$Q$-polynomial of $G$ is
\[
f_Q (\lambda ) = (f_{Q_2 } (\lambda  -2))^{m_1 } f_{Q_1 } \left(( -
1)^{n_2}\frac{(\lambda -r_1 n_2)\cdot {f_{Q_2} (\lambda )}}{
{f_{\overline Q_2 } ( {n_2 - \lambda } )}} \right)\left( ( -
1)^{n_2} \frac{{f_{\overline Q_2 } \left( {n_2 - \lambda }
\right)}}{{f_{Q_2} (\lambda )}}\right)^{n_1}.
\]}
\\

The following Proposition 2.13 exhibits the $Q$-coronal of the join
$G_1\vee G_2$ of two regular graphs $G_1$ and $G_2$.\\
\\
\textbf{Proposition 2.13.} {\it Let $G_1$ be an $r_1$-regular graph
on $n_1$ vertices and $G_2$ be an $r_2$-regular graph on $n_2$
vertices. Also let $G=G_1\vee G_2$. Then
\[
\Gamma _Q (\lambda ) = \frac{{(\lambda  - n_2  - 2r_1 )n_2  +
(\lambda  - n_1  - 2r_2 )n_1  + 2n_1 n_2 }}{{(\lambda  - n_2  - 2r_1
)(\lambda  - n_1  - 2r_2 ) - n_1 n_2 }}.
\]}\\
\textbf{Proof.} This follows directly from Theorem 2.2, Corollary
2.4 and Proposition 2.10. Namely, from Theorem 2.2, we have
\[
f_{\overline Q_i } (\lambda ) = ( - 1)^{n_i} \left( {1 +
\frac{n_i}{{n_i - 2 - 2r_i - \lambda }}} \right)f_{Q_i} (n_i - 2 -
\lambda ),\quad (i=1,2).
\]
It follows from $f_{\overline Q } (\lambda ) =f_{\overline Q_1 }
(\lambda )f_{\overline Q_2 } (\lambda )$ that
\[
f_{\overline Q } (n - 2 - \lambda ) = ( - 1)^{n_1  + n_2 } \left( {1
+ \frac{{n_1 }}{{\lambda  - n_2  - 2r_1 }}} \right)\left( {1 +
\frac{{n_2 }}{{\lambda  - n_1  - 2r_2 }}} \right)f_{Q_1 } (\lambda -
n_2 )f_{Q_2 } (\lambda  - n_1 ).
\]
By Corollary 2.4, one gets
\[
f_{Q}(\lambda)  = \left({1 - \frac{{n_1 n_2 }}{{(\lambda  - n_1  -
2r_2 )(\lambda  - n_2  - 2r_1 )}}} \right)f_{Q_1 } (\lambda - n_2
)f_{Q_2 } (\lambda  - n_1 ).
\]
Now the result follows easily from Proposition 2.10. $\Box$\\

Next, we shall derive the $Q$-generating function for some graphs
obtained by the use of some operation on graphs, such as the
complement of a graph, the direct sum and the join of two
graphs and so on.\\
\\
\textbf{Theorem 2.14.} {\it For the generating function $W_Q(t)$ for
the numbers $N_k$ of semi-edge walks of length $k$ in a graph $G$,
we have
\begin{equation}\label{5}
W_{\overline Q } (t) = \frac{{ - W_Q \left( {\frac{t}{{(n - 2)t -
1}}} \right)}}{{(n - 2)t - 1 + tW_Q \left( {\frac{t}{{(n - 2)t -
1}}} \right)}},
\end{equation}
\begin{equation}\label{6}
W_{Q_1  \oplus Q_2 } (t) = W_{Q_1 } (t) + W_{Q_2 } (t),
\end{equation}
\begin{equation}\label{7}
W_{Q_1  \vee Q_2 } (t) = \frac{M}{{1 - tM}},
\end{equation}
where
\[
M = \sum\limits_{i = 1}^2 {\frac{{W_{Q_i } (t)}}{{(n_i  - n)t + 1 +
tW_{Q_i } (t)}}}.
\]}\\
\\
\textbf{Proof.} From Proposition 2.1, one has
\[
W_{\overline Q } (t) = \frac{1}{t}\left( {( - 1)^n \frac{{f_Q \left(
{n - 2 - \frac{1}{t}} \right)}}{{f_{\overline Q } \left(
{\frac{1}{t}} \right)}} - 1} \right)
\]
and
\[
W_Q \left( {\frac{1}{{n - 2 - \frac{1}{t}}}} \right) = \left( {n - 2
- \frac{1}{t}} \right)\left( {( - 1)^n \frac{{f_{\overline Q }
\left( {\frac{1}{t}} \right)}}{{f_Q \left( {n - 2 - \frac{1}{t}}
\right)}} - 1} \right),
\]
which implies that the required result (\ref{5}). The formula
(\ref{6}) is obvious. Next, we shall prove (\ref{7}). According to
(\ref{5}) and (\ref{6}), one gets
\begin{align}\label{8}
W_{Q_1  \vee Q_2 } (t) &= W_{\overline {\overline {Q_1 }  \oplus
\overline {Q_2 } } } (t)\nonumber\\
& = \frac{{ - W_{\overline {Q_1 }  \oplus \overline {Q_2 } } \left(
{\frac{t}{{(n - 2)t - 1}}} \right)}}{{(n - 2)t - 1 + tW_{\overline
{Q_1 }  \oplus \overline {Q_2 } } \left( {\frac{t}{{(n - 2)t - 1}}}
\right)}}\nonumber\\
& =  - \frac{{W_{\overline {Q_1 } } \left( {\frac{t}{{(n - 2)t -
1}}} \right) + W_{\overline {Q_2 } } \left( {\frac{t}{{(n - 2)t -
1}}} \right)}}{{(n - 2)t - 1 + t\left( {W_{\overline {Q_1 } } \left(
{\frac{t}{{(n - 2)t - 1}}} \right) + W_{\overline {Q_2 } } \left(
{\frac{t}{{(n - 2)t - 1}}} \right)} \right)}}.
\end{align}
For $i=1,2$, the formula (\ref{5}) implies that
\begin{align}\label{9}
W_{\overline {Q_i } } \left( {\frac{t}{{(n - 2)t - 1}}} \right) &=
\frac{{ - W_{Q_i } \left( t \right)}}{{(n_i  - 2)\frac{t}{{(n - 2)t
- 1}} - 1 + \frac{t}{{(n - 2)t - 1}}W_{Q_i } \left( t
\right)}}\nonumber\\& = \frac{{ - ((n - 2)t - 1)W_{Q_i } \left( t
\right)}}{{(n_i  - n)t + 1 + tW_{Q_i } \left( t \right)}}.
\end{align}
Substituting (\ref{9}) back into (\ref{8}), we obtain the
required result (\ref{7}). $\Box$\\
\\
\textbf{Remark 2.15.}  In view of the formulas (\ref{5}) and
(\ref{6}), the formulas (\ref{7}) may generalized  to the case
$k>2$, that is,
\begin{equation}\label{10}
W_{Q_1  \vee Q_2 \vee\cdots\vee Q_k} (t) = \frac{M}{{1 - tM}},
\end{equation}
where
\[
M = \sum\limits_{i = 1}^k {\frac{{W_{Q_i } (t)}}{{(n_i  - n)t + 1 +
tW_{Q_i } (t)}}}.
\]\\
\textbf{Example 2.16.} Consider the complete multipartite graph
$G=K_{n_1,n_2,\ldots,n_k}$, which can be represented as the join of
graphs $G_1, G_2\cdots, G_k$, all of which contain only isolated
vertices. For an $r$-regular graph $G$, its $Q$-generating function
$W_Q(t)=\frac{n}{1-2rt}$ for the numbers $N_k$ of semi-edge walks of
length $k$ in $G$ (see the proof of Theorem 2.2). Hence the
$Q_i$-generating function for the numbers $N_k$ of semi-edge walks
of length $k$ in $G_i$ equals $n_i$ for $i=1,2,\ldots,k$. Now
applying (\ref{10}), we obtain
\[
W_Q (t) = \left( {\left( {\sum\limits_{i = 1}^k {\frac{{n_i }}{{(n_i
- n)t + 1 + tn_i }}} } \right)^{ - 1}  - t} \right)^{ - 1}.
\]
According to the proof of Proposition 2.10, the $Q$-coronal of the
complete multipartite graph $G=K_{n_1,n_2,\ldots,n_k}$ equals
\[
\Gamma _Q (\lambda ) = \left( {\left( {\sum\limits_{i = 1}^k
{\frac{{n_i }}{{\lambda  - n + 2n_i }}} } \right)^{ - 1}  - 1}
\right)^{ - 1}.
\]\\

Finally, Liu and Lu \cite{Liu2013} introduced the definitions of the
subdivision-vertex neighbourhood corona and subdivision-edge
neighbourhood corona for two graphs $G_1$ and $G_2$, and their
$Q$-polynomials are determined by using the $Q$-coronal of $G_2$.
Clearly, Applying the combinatorial interpretation of the
$Q$-coronal of graphs (see Proposition 2.10), we may obtain many
alternative calculations of the $Q$-polynomials of the
subdivision-vertex neighbourhood corona and subdivision-edge
neighbourhood corona of $G_1$ and $G_2$. These contents are
omitted.\\
\\
\textbf{Acknowledgements} This work was partially supported by the
National Natural Science Foundation of China (No. 11271334), the
Natural Science Foundation of Zhejiang Province, China (No.
LY12A01006) and the Scientific Research Fund of Zhejiang Provincial
Education Department (No. Y201225862).

\end{document}